\documentclass{amsart}
\usepackage{amssymb}
\copyrightinfo{2013}{V.X. Genest \emph{et al.}}

\catcode`,\active

\catcode`\,12

\catcode`,\active

\catcode`\,12

\newcommand*\pGqskip{8mu}
\catcode`,\active
\newcommand*\pGq{\begingroup
        \catcode`\,\active
        \def ,{\mskip\pGqskip\relax}%
        \dopGq
}
\catcode`\,12
\def\dopGq#1#2#3#4#5{%
        {}_{#1}F_{#2}\bigg(\genfrac..{0pt}{}{#3}{#4}\,;\,#5\biggr)%
        \endgroup
}
\newcommand{\pd}{\partial}

\theoremstyle{definition}

\theoremstyle{remark}
\newtheorem*{remark}{Remark}

\numberwithin{equation}{section}
\begin{document}
\title[Tridiagonalization and the Racah--Wilson algebra]{Tridiagonalization of the hypergeometric operator and the Racah--Wilson algebra}
\author[V.X. Genest]{Vincent X. genest}
\address{Centre de Recherches Math\'ematiques, Universit\'e de Montr\'eal, P.O. Box 6128, Centre-ville Station, Montr\'eal (Qu\'ebec) H3C 3J7}
\email{genestvi@crm.umontreal.ca}
\author[M. E. H. Ismail]{Mourad E. H. Ismail}
\address{Department of Mathematics, King Saud University, Riyadh, Saudi Arabia and Department of Mathematics, University of Central Florida, Orlando, Florida 32816, USA}
\email{mourad.eh.ismail@gmail.com}
\author[L. Vinet]{Luc Vinet}
\address{Centre de Recherches Math\'ematiques, Universit\'e de Montr\'eal, P.O. Box 6128, Centre-ville Station, Montr\'eal (Qu\'ebec) H3C 3J7}
\email{vinetl@crm.umontreal.ca}
\author[A. Zhedanov]{Alexei Zhedanov}
\address{Centre de Recherches Math\'ematiques, Universit\'e de Montr\'eal, P.O. Box 6128, Centre-ville Station, Montr\'eal (Qu\'ebec) H3C 3J7}
\email{zhedanov@yahoo.com}
\date{}
\dedicatory{}
\begin{abstract}
The algebraic underpinning of the tridiagonalization procedure is investigated. The focus is put on the tridiagonalization of the hypergeometric operator and its associated quadratic Jacobi algebra. It is shown that under tridiagonalization, the quadratic Jacobi algebra becomes the quadratic Racah--Wilson algebra associated to the generic Racah/Wilson polynomials. A degenerate case leading to the Hahn algebra is also discussed.
\end{abstract}
\maketitle
\section{Introduction}
This paper aims to offer an algebraic framework for the tridiagonalization approach to the study of orthogonal polynomials. This procedure has been investigated by Ismail and Koelink in a series of papers \cite{2011_Ismail&Koelink_AdvApplMath_46_379,2012_Ismail&Koelink_SIGMA_8_61,2012_Ismail&Koelink_AnalAppl_10_327}. In \cite[Section 3]{2012_Ismail&Koelink_AnalAppl_10_327}, these authors have shown in particular that a special case of the Wilson polynomials can be obtained from the Jacobi polynomials by tridiagonalization. The hyper\-geo\-metric polynomials of the Askey tableau can be described from the most general to the more specific by starting with the 4-parameter Wilson polynomials, which sit atop the tableau \cite[Chapter 9]{2010_Koekoek_&Lesky&Swarttouw}. It is however appealing to develop an understanding of the Askey scheme from the bottom up by constructing the more involved families of polynomials from the simpler ones. This is in part what the tridiagonalization method allows to do. The idea of ``lifting'' within the Askey scheme is also in that spirit \cite{2010_Atakishiyeva&Atakakishiyev_JPhysA_43_145201}, \cite{1996_Berg&Ismail_CanJMath_48_43} \cite[Sections 15.1, 15.2]{Ismail-2009}.

The main concept behind the tridiagonalization approach is to start from an operator $L$ that has a family of orthogonal polynomials $\{P_{n}(x)\}$ as eigenfunctions and to combine it with other operators so as to obtain a new operator $M$ that acts on $\{P_{n}(x)\}$ in a tridiagonal fashion. The diagonalization of the operator $M$ by another family of orthogonal polynomials $\{Q_{n}(x)\}$, when possible, will relate $\{P_{n}(x)\}$ and $\{Q_{n}(x)\}$, whose characterization will benefit from the setup.

It is known that many of the important properties of the classical orthogonal polynomials are encoded in polynomial algebras and their representations \cite{1991_Zhedanov_TheorMathPhys_89_1146}. This is done for instance by using as generators the two operators involved in the bispectrality of the family of polynomials of interest.  We here integrate the tridiagonalization method into this algebraic picture. We focus on the connection between  Jacobi and Wilson polynomials via tridiagonalization and explain how the Wilson/Racah algebra emerges in the process.

The outline of the paper is as follows. The tridiagonalization of the hypergeometric operator $L$ is revisited in Section 2. The approach of \cite{2012_Ismail&Koelink_AnalAppl_10_327} is generalized by considering both left and right multiplication by $x$. It is observed that $L$ and the tridiagonalized operator $M$ interchange their role under the transformation $x\rightarrow x^{-1}$. In Section 3, it is shown that the generic Wilson polynomials with 4 parameters arise from the tridiagonalization of $L$. This leads naturally to Koornwinder's integral representation for the Wilson polynomials \cite{1985_Koornwinder_OPs&Appl_174}. The case when both $L$ and $M$ preserve the space of polynomials of degree smaller or equal to a positive integer $N$ is considered. This finite-dimensional reduction is the object of Section 4. The connection with the Wilson/Racah algebra is explained in Section 5, where it is shown that the hypergeometric generator $L$ and the tridiagonalized operator $M$ realize this algebra. Reference is made to the realization of the Wilson/Racah algebra out of the intermediate Casimir operators occurring in the coupling of three irreducible representation of $\mathfrak{su}(1,1)$. The operators introduced in the tridiagonalization of $L$ will be observed to coincide with the Casimir operators arising in the one-dimensional model of the Racah problem elaborated in \cite{2014_Genest&Vinet&Zhedanov_JPhysA_47_025203}. Section 6 is dedicated to the degenerate case where the parameters of the tridiagonalized operator $M$ are such that it becomes a first order differential operator. It is then seen that the Hahn algebra and polynomials occur.
\section{Tridiagonalization procedure for the hypergeometric operator}
In this section, the tridiagonalization procedure for the hypergeometric operator is revisited and generalized. The quadratic Jacobi algebra associated to this operator is also presented.
\subsection{The classical Jacobi polynomials and the quadratic Jacobi algebra}
Consider the ordinary hypergeometric operator 
\begin{equation}
\label{hyp_op}
L= x(1-x) \partial_x^2 + \big(\alpha+1 - (\alpha+\beta+2)\,x\big) \partial_x,
\end{equation}
where $\alpha, \beta$ are real parameters satisfying $\alpha,\beta >-1$. The operator $L$ possesses an infinite family of polynomial eigenfunctions $\{\widehat{P}_{n}^{(\alpha,\beta)}(x)\}_{n=0}^{\infty}$ that obey the eigenvalue equation
\begin{equation}
 \label{L_Jac}
L\,\widehat{P}_n^{(\alpha, \beta)}(x) = \lambda_n \, \widehat{P}_n^{(\alpha, \beta)}(x),
\end{equation}
with eigenvalues
\begin{equation}
\label{lambda_Jac}
\lambda_n = -n(n+\alpha+\beta+1).
\end{equation}
The monic polynomials
\begin{equation}
\label{Jac_hyp} 
\widehat{P}_n^{(\alpha, \beta)}(x) =\frac{(-1)^n(\alpha+1)_n}{(\alpha+\beta+n+1)_n}\; \pGq{2}{1}{-n,n+\alpha+\beta+1}{\alpha+1}{x}, 
\end{equation}
with $(a)_n$ denoting the Pochhammer symbol
\begin{equation*}
 (a)_{n}=a(a+1)(a+2)\cdots(a+n-1),\qquad (a)_0=1,
\end{equation*}
coincide, up to an affine transformation, with the Jacobi polynomials \cite{2001_Andrews&Askey&Roy}. The polynomials $\widehat{P}_n^{(\alpha, \beta)}(x) $ satisfy the three-term recurrence relation
\begin{equation}
\label{3_t_Jac}
\widehat{P}_{n+1}^{(\alpha, \beta)}(x) + b_n \widehat{P}_n^{(\alpha, \beta)}(x) + u_n \widehat{P}_{n-1}^{(\alpha, \beta)}(x) = x \widehat{P}_n^{(\alpha, \beta)}(x), 
\end{equation}
with the recurrence coefficients 
\begin{align*}
u_n&=\frac{n(n+\alpha)(n+\beta)(n+\alpha+\beta)}{(2n+\alpha+\beta-1)(2n+\alpha+\beta)^2(2n+\alpha+\beta+1)}, 
\\
b_n&= \frac{1}{2} + \frac{\alpha^2-\beta^2}{4} \left(\frac{1}{2n+\alpha+\beta} -  \frac{1}{2n+\alpha+\beta+2} \right). 
\end{align*}
They obey the orthogonality relation
\begin{equation}
\label{ort_Jac}
\int_0^1  \widehat{P}_n^{(\alpha, \beta)}(x) \widehat{P}_m^{(\alpha, \beta)}(x)\;x^{\alpha} (1-x)^{\beta} dx = h_n \: \delta_{nm},
\end{equation}
where the normalization factor $h_n$ is given by
\begin{equation}
\label{h_n_Jac}
h_n = \frac{\Gamma(\alpha+1) \Gamma(\beta+1)}{\Gamma(\alpha+\beta+2)} u_1 u_2 \dots u_n,
\end{equation}
where $\Gamma(x)$ is the Gamma function.

There is an associative algebra known as the \emph{quadratic Jacobi algebra} which is associated to the hypergeometric operator and to the Jacobi polynomials \cite[Formula $(3.4)$ and Section 4]{1992_Granovskii&Lutzenko&Zhedanov_AnnPhys_217_1}. This algebra can be obtained as follows. Let $K_1$ and $K_2$ be defined as 
\begin{align}
\label{Rea-1}
 K_1=L,\qquad K_2=X,
\end{align}
where $L$ is given by \eqref{hyp_op} and where $X$ is the ``multiplication by $x$'' operator; that is $X\,f(x)= x f(x)$. Introduce a third operator $K_3$, defined as $K_3\equiv [K_1,K_2] $, where $[a,b]=ab-ba$ is the commutator. From \eqref{hyp_op}, it follows that $K_3$ has the expression
\begin{align*}
 K_3=-2x(x-1)\pd_{x}+(1+\alpha-(\alpha+\beta+2)x).
\end{align*}
By a direct calculation, one finds that $K_1$, $K_2$ and $K_3$ satisfy the relations
\begin{align}
\label{QR3}
 \begin{aligned}
  &[K_1,K_2]=K_3,
  \\
  &[K_2,K_3]=a_2 K_2^2+d K_2,
  \\
  &[K_3,K_1]=a_2 \{K_1,K_2\}+d K_1+c_2 K_2+e_2,
 \end{aligned}
\end{align}
where $\{a,b\}=ab+ba$ and where the structures constants $a_2,d, c_2, e_2$ are given by
\begin{align*}
 a_2=2, \quad d=-2,\quad c_2=-(\alpha+\beta)(\alpha+\beta+2),\quad e_2=(\alpha+1)(\alpha+\beta).
\end{align*}
The relations \eqref{QR3} define the Jacobi algebra. This algebra has a Casimir operator
\begin{align}
\label{QR3-Cas}
 Q=a_2 \{K_2^2, K_1\}+K_3^2+(a_2^2+c_2) K_2^2+d \{K_1,K_2\}+(d a_2+2 e_2)K_2,
\end{align}
which commutes with all generators $K_1, K_2, K_3$. In the realization \eqref{Rea-1}, one has
\begin{align}
\label{value}
 Q\,f(x)=(\alpha^2-1) f(x).
\end{align}
The polynomials $\{\widehat{P}_{n}^{(\alpha,\beta)}(x)\}_{n=0}^{\infty}$ form a basis for infinite-dimen\-sional irreducible representations of the quadratic Jacobi algebra \eqref{QR3}. In this basis, the operator $K_1$ is diagonal with eigenvalues \eqref{lambda_Jac}, and the operator $K_2$ is tridiagonal, as is manifest from \eqref{3_t_Jac}. The monomials $\{x^{n}\}_{n=0}^{\infty}$ also form a basis for these irreducible representations. In this basis, both generators $K_1$ and $K_2$ are two-diagonal. Indeed, one has
\begin{equation}
\label{L_2_diag}
K_1\,x^n =L\,x^{n}= \lambda_n x^n + \epsilon_n x^{n-1},\qquad \text{with}\qquad \epsilon_n =n(n+\alpha).
\end{equation}
and evidently $X\,x^{n}=x^{n+1}$.
\subsection{Tridiagonalization of $L$}
Consider the operator
\begin{equation}
 \label{M-Def}
M = \tau_1 X L + \tau_2 L X + \tau_3 X + \tau_0,
\end{equation}
where $\tau_i, \: i=0,1,2,3$ are arbitrary parameters such that $\tau_1+\tau_2 \ne 0$. Without loss of generality one can take $\tau_1+\tau_2=1$, as this corresponds to scaling the operator $M$. Upon using \eqref{hyp_op} in \eqref{M-Def}, one finds that $M$ has the expression
\begin{multline}
 \label{expl_M}
 M = x^2(1-x)\;\pd_x^2 + x \left(\alpha+1 +2\,\tau_2 -(\alpha+\beta+2+2\tau_2)x   \right) \partial_x
 \\
 -(\tau_2 (\alpha+\beta+2)-\tau_3)\,x +(\alpha+1)\tau_2 +\tau_0.
\end{multline}
\begin{remark}
The condition $\tau_1 + \tau_2 \ne 0$ is essential for the operator $M$ to be a second-order differential operator. If $\tau_1+\tau_2=0$, $M$ becomes a first-order operator.
\end{remark}
From \eqref{L_Jac}, \eqref{3_t_Jac}, \eqref{L_2_diag} and \eqref{M-Def}, it is seen that the operator $M$ has the following properties:
\begin{itemize}
 \item It acts in a three-diagonal fashion on the polynomials $\widehat{P}_n^{(\alpha,\beta)}(x)$. One has indeed
\begin{multline}
\label{M_P_Jac}
 M \widehat{P}_n^{(\alpha,\beta)}(x) = (\tau_1 \lambda_{n} + \tau_2 \lambda_{n+1} + \tau_3) \widehat{P}_{n+1}^{(\alpha,\beta)}(x)    \\
 +(\lambda_nb_n + \tau_3\,b_n+ \tau_0 )\widehat{P}_n^{(\alpha,\beta)}(x)
 + (\tau_1 \lambda_{n} + \tau_2 \lambda_{n-1} + \tau_3)\,u_n \,\widehat{P}_{n-1}^{(\alpha,\beta)}(x).
\end{multline}
\item It acts in a two-diagonal fashion on the monomials $x^n$. One has
\begin{equation}
 \label{M_x_n}
M x^n = (\tau_1 \lambda_n + \tau_2 \lambda_{n+1} + \tau_3)\,x^{n+1} + (\tau_1 \epsilon_n + \tau_2 \epsilon_{n+1} + \tau_0)\,x^n,
\end{equation}
where $\lambda_n$ and $\epsilon_n$ are respectively given by \eqref{lambda_Jac} and \eqref{L_2_diag}.
\end{itemize}
There is an interesting duality between the operators $L$ and $M$. Indeed, one observes using the explicit expressions \eqref{L_Jac} and \eqref{M-Def} with the condition $\tau_1+\tau_2=1$ that the operator $L$ can be expressed in terms of the operator $M$ as follows:
\begin{equation}
 \label{rec_LM}
L= \tau_1 X^{-1} M + \tau_2 M X^{-1} -(\tau_0 + 2 \tau_1 \tau_2) X^{-1} + 2 \tau_1 \tau_2 - \tau_3,
\end{equation}   
where $X^{-1}$ is the ``multiplication by $x^{-1}$'' operator, i.e. $X^{-1}\,f(x)=x^{-1}\,f(x)$. The relations \eqref{M-Def} and \eqref{rec_LM} have the same structure under the substitution $X \to X^{-1}$. Consequently, one can expect that the operator $L$ will be three-diagonal in an appropriate (discrete) basis of eigenfunctions of $M$.  It is convenient to introduce two new parameters $\nu_1$ and $\nu_2$ instead of $\tau_1$ and $\tau_3$
\begin{equation}
\label{tau_nu}
\tau_2= \frac{1}{2} (1+\nu_1-\nu_2), \qquad \tau_3 = \nu_1\nu_2 +\frac{(\alpha+\beta)(\nu_1+\nu_2-1)}{2}.
\end{equation}
Consider the change of variable $x = 1/y$ together with the similarity transformation
\begin{align}
\label{Sim}
 \widetilde{L}=y^{\nu_2-1}\,L\,y^{1-\nu_2},\qquad \widetilde{M}=y^{\nu_2-1}\,M\,y^{1-\nu_2}.
\end{align}
Under these transformations, the operator $L$ takes the form
\begin{equation}
\label{t_L}
\widetilde{L}=y^2 (y-1) \partial_y^2 + y(a_1y + b_1)\partial_y + c_1y + d_1,
\end{equation}
where the coefficients $a_1,b_1, c_1, d_1$ are given by
\begin{alignat*}{2}
 &a_1=3-\alpha-2\nu_2, &\qquad b_1&=\alpha+\beta+2\nu_2-2,\quad 
 \\
 &c_1 =(\nu_2-1)(\alpha+\nu_2-1), &\qquad d_1&=-(\nu_2-1)(\alpha+\beta+\nu_2).
\end{alignat*}
The operator $M$ becomes
\begin{equation}
\label{t_M}
\widetilde{M}= y(y-1) \partial_y^2 + (a_2y + b_2) \partial_y + d_2, 
\end{equation}
where $a_2, b_2, d_2$ read
\begin{gather*}
 a_2= 2-\nu_1-\nu_2-\alpha, \qquad b_2=\alpha+\beta+\nu_1+\nu_2-1, 
 \\
 d_2=\tau_0 + \nu_1 \nu_2 +\frac{1}{2}(\alpha-1) (\nu_1+\nu_2-1).
\end{gather*}
As is seen from \eqref{t_L} and \eqref{t_M}, the operators $L$ and $M$ have their roles interchanged under the change of variable $y=1/x$ and the similarity transformation \eqref{Sim}. Indeed the operator $\widetilde{M}$ coincides with the hypergeometric operator (up to an affine transformation), while $\widetilde{L}$ has a form similar to $M$.  This means in particular that there exists a set of eigenfunctions $\{\psi_{n}(x)\}$ for $n=0,1,2,\ldots$ that satisfy the eigenvalue equation
\begin{equation*}
M\,\psi_n(x) = \widetilde{\lambda}_n \psi_n(x) ,
\end{equation*}
and on which the operator $L$ acts in a three-diagonal fashion
\begin{equation}
\label{L_psi}
L\,\psi_n(x) = \xi_n \psi_{n+1}(x) + \eta_n \psi_n(x) + \zeta_n \psi_{n-1}(x) ,
\end{equation}
for some coefficients $\xi_n, \eta_n, \zeta_n$. Indeed, since $\widetilde{M}$ is the ordinary hypergeometric operator, it possesses a set of polynomials as eigenfunctions. It is verified that $M$ admits the discrete set of eigenfunctions
\begin{equation}
\label{psi_M}
\psi_n(x) = x^{\nu_2-1} \widehat{P}_n^{(\widetilde{\alpha}, \widetilde{\beta)}}(1/x), \quad \widetilde{\lambda}_n = n(n-1) + a_2 n + d_2, 
\end{equation}   
where $\widetilde{\beta}=\beta, \; \widetilde{\alpha}= -\alpha-\beta-\nu_1 - \nu_2$. The property \eqref{L_psi} then follows from the fact that $\widehat{P}_n^{(\alpha,\beta)}(x)$ are the polynomials \eqref{Jac_hyp} satisfying the three-term recurrence relation \eqref{3_t_Jac}.
\section{Wilson polynomials arising from tridiagonaliation}
In this section, we show that the tridiagonalization of the hypergeometric operator $L$ leading to the operator $M$ gives rise to the generic Wilson polynomials.

Let $\psi(x; \Lambda)$ be an eigenfunction of $M$ corresponding to the eigenvalue $\Lambda$
\begin{equation}
 \label{M_psi}
M \psi(x;\Lambda) = \Lambda\,\psi(x;\Lambda),
\end{equation}
which is integrable with respect to the Jacobi weight. Consider the expansion of $\psi(x;\Lambda)$ in the formal Jacobi series
\begin{equation}
\label{psi_P}
\psi(x;\Lambda) = \sum_{k=0}^{\infty} G_k(\Lambda) \widehat{P}_k^{(\alpha,\beta)}(x).
\end{equation}
Let us write the expansion coefficients as
\begin{equation*}
G_k(\Lambda) = G_0(\Lambda)\,\Xi_k\,Q_k(\Lambda), 
\end{equation*}
where $\Xi_k$ are normalization coefficients, $Q_{0}(\Lambda)=1$ and where  $G_0(\Lambda)$ is given by
\begin{equation*}
G_0(\Lambda) = \frac{1}{h_0} \: \int_0^1 \psi(x;\Lambda) x^{\alpha} (1-x)^{\beta} dx.
\end{equation*}
Upon applying $M$ on both sides of \eqref{psi_P} and using \eqref{M_P_Jac}, one finds that the coefficients $Q_k(\Lambda)$ satisfy the three-term recurrence relation
\begin{equation*}
E_n^{(1)} Q_{n+1}(\Lambda) + E_n^{(2)} Q_n(\Lambda) + E_n^{(3)} Q_{n-1}(\Lambda) = \Lambda \,Q_n(\Lambda), 
\end{equation*}
where
\begin{gather*}
 \begin{gathered}
E_n^{(1)}=u_{n+1} \frac{\Xi_{n+1}}{\Xi_n} (\tau_1 \lambda_{n+1} + \tau_2 \lambda_n + \tau_3),\qquad 
E_n^{(2)}=\lambda_n b_n + \tau_3 b_n + \tau_0, 
\\
E_n^{(3)}=\frac{\Xi_{n-1}}{\Xi_n} (\tau_1 \lambda_{n-1} + \tau_2 \lambda_n + \tau_3).
\end{gathered}
\end{gather*}
Let us choose the coefficients $\Xi_n$ to ensure that
\begin{equation*}
E_n^{(1)}= u_{n+1} \frac{\Xi_{n+1}}{\Xi_n} (\tau_1 \lambda_{n+1} + \tau_2 \lambda_n + \tau_3) =1,
\end{equation*}
and that $\Xi_0=1$. These requirements determine $\Xi_n$ uniquely. One then finds that $Q_n(\Lambda)$ are polynomials in $\Lambda$ with the three-term recurrence relation 
\begin{equation}
\label{rec_Q} 
Q_{n+1}(\Lambda) + B_n\,Q_n(\Lambda) + U_n\,Q_{n-1}(\Lambda) = \Lambda\,Q_n(\Lambda),
\end{equation}
recurrence coefficients
\begin{equation*}
B_n = E_n^{(2)}=   \lambda_n b_n + \tau_3 b_n + \tau_0 ,\quad U_n = u_n  (\tau_1 \lambda_{n-1} + \tau_2 \lambda_n + \tau_3)(\tau_1 \lambda_{n} + \tau_2 \lambda_{n-1} + \tau_3),
\end{equation*}
and initial conditions $Q_{-1}(\Lambda)=0$ and $Q_{0}(\Lambda)=1$. It follows that $Q_n(\Lambda)$ are orthogonal polynomials in $\Lambda$. These polynomials can be identified with the Wilson polynomials $W_n(x;a_1,a_2,a_3,a_4)$, which depend on 4 parameters $a_1, a_1, a_3, a_4$ and which are defined by the three-term relation
\begin{equation}
 \label{rec_W}
W_{n+1}(x) + (A_n +C_n - a_1^2)\,W_n(x) + A_{n-1} C_n W_{n-1}(x) = xW_n(x), \end{equation}
where
\begin{align*}
 \begin{aligned}
  A_n&=\frac{(n+g-1)(n+a_1+a_2)(n+a_1+a_3)(n+a_1+a_4)}{(2n +g-1)(2n+g)},
  \\
  C_n&=\frac{n(n+a_2+a_3-1)(n+a_2+a_4-1)(n+a_3+a_4-1)}{(2n+g-2)(2n +g-1)},
 \end{aligned}
\end{align*}
and where  $g=a_1+a_2+a_3+a_4$. Upon comparing \eqref{rec_Q} and \eqref{rec_W}, one finds 
\begin{equation}
\label{Q_W_id}
Q_n(\Lambda) = W_n(\gamma-\Lambda;a_1,a_2,a_3,a_4),
\end{equation}
with $\gamma = \frac{1}{2} \left( a_1 + a_2 - a_1^2 - a_2^2 +2\tau_0\right)$ and where
\begin{alignat*}{2}
\alpha&=a_1+a_2-1, &\qquad \beta&=a_3+a_4-1,
\\
\nu_1&= 1-a_1-a_3, &\qquad \nu_2&=1-a_2-a_3.
\end{alignat*}
We have thus identified the polynomials $Q_n(\Lambda)$ with the generic Wilson polynomials. This indicates that the tridiagonalization procedure, when applied to the hypergeometric operator $L$, leads to the generic Wilson polynomials, which depend on four parameters. 

\begin{remark}
While four parameters appear in the definition \eqref{M-Def} of $M$, only two of them are essential since $\tau_0$ is merely a translation and since $M$ can be scaled by an arbitrary factor. Together with the two parameters $\alpha$, $\beta$ of the hypergeometric operator, this accounts for the four arbitrary parameters of the Wilson polynomials arising in \eqref{Q_W_id}.
\end{remark}

Let us parametrize the eigenvalues $\Lambda$ appearing in \eqref{M_psi} as
\begin{equation*}
\Lambda= \tau_0 + \frac{1-\alpha}{2} + (\alpha+q)(q+1) + \frac{(\nu_1-\nu_2)(\alpha+2q+1)}{2} \label{Lambda_q}, \end{equation*}
where $q$ is a new parameter. It is verified that \eqref{M_psi} admits the solution
\begin{equation}
\label{Dompe}
\psi(x;\Lambda) = x^{q} (1-x)^{-\beta}\;\pGq{2}{1}{\alpha_1,\beta_1}{\gamma_1}{x},
\end{equation}
where $\alpha_1$, $\beta_1$ and $\gamma_1$ are given by
\begin{equation*}
\alpha_1 =\alpha+1 +\nu_1+q, \qquad \beta_1 = 1-\beta-\nu_2+q, \qquad \gamma_1 =  2 +2q +\alpha + \nu_1 - \nu_2 \label{al_1}. \end{equation*}
Note that one has $\gamma_1-\alpha_1-\beta_1=\beta$ and hence it follows that for $q>0$ and $\beta>0$ the integrals
\begin{equation*}
\int_0^1 x^{\alpha} (1-x)^{\beta} x^n \psi(x;\Lambda) dx,
\end{equation*}
exist for all non-negative integer $n$. Thus the function $\psi(x;\Lambda)$ taken as in \eqref{Dompe} satisfies the required property of being integrable with respect to the Jacobi weight function. The other independent solution of \eqref{M_psi} does not have this property and is hence discarded.  Formally using \eqref{psi_P} suggests the following integral representation of the generic Wilson polynomials $Q_n(\Lambda)$:
\begin{equation*}
Q_n(\Lambda) = \frac{1}{h_n G_0(\Lambda) \Xi_n} \: \int_0^1 \psi(x;\Lambda) \widehat{P}_k^{(\alpha,\beta)}(x) x^{\alpha} (1-x)^{\beta} dx. 
\end{equation*}
This formula is seen to coincide with Koornwinder's integral representation for the Wilson polynomials \cite{1985_Koornwinder_OPs&Appl_174}. It is interesting to note that the dual solutions $\psi_n(x)$ described by  \eqref{psi_M} do not satisfy the condition that the integrals
\begin{equation*}
\int_0^1 \psi_n(x) x^{\alpha} (1-x)^{\beta} dx
\end{equation*}
exist for all non negative integers $n$. Indeed, the function $\psi_n(x)$ has the polynomial part $\widehat{P}_n^{(\widetilde{\alpha},\widetilde{\beta})}(1/x)$. Hence, for sufficiently large $n$, it is seen that the integral 
\begin{equation*}
\int_0^1 x^{-n} x^{\alpha} (1-x)^{\beta} dx
\end{equation*}
will diverge.
\begin{remark}
In the process of describing the algebraic underpinning of the tridiagonalization procedure, we have seen the integral representation of the Wilson polynomials of \cite{1985_Koornwinder_OPs&Appl_174} arise as a side result. As pointed out by a referee, one could have started instead from that integral representation to derive the embedding of the Racah--Wilson algebra into the Jacobi algebra.
\end{remark}

\section{Finite-dimensional reduction}
In general, the two discrete sets of eigenfunctions $\{\widehat{P}_n^{(\alpha,\beta)}(x)\}$ and $\{\psi_n(x)\}$ of the operators $L$ and $M$ are infinite-dimensional. However, there is a case for which it is possible to work with finite-dimensional sets. This is the object of this section.

Consider the action of the operators $L$ and $M$ on the monomials $x^n$. In view of \eqref{L_2_diag}, it is manifest that $L$ preserves the space of polynomials of a given degree $N$, where $N$ is a non-negative integer. In contrast, it is visible from \eqref{M_x_n} that in general $M$ does not preserve the space of polynomials of a given degree, as it maps polynomials of degree $n$ to polynomials of degree $n+1$. Consider the special case for which the condition 
\begin{equation}
\tau_1 \lambda_N + \tau_2 \lambda_{N+1} + \tau_3 =  \tau_1 (\lambda_N -\lambda_{N+1}) +\lambda_{N+1} + \tau_3  =0, \label{trunc_M}
\end{equation}
holds for some non-negative integer $N$. It is clear from \eqref{M_x_n} that in this case the operator $M$ preserves the space of polynomials of degree less or equal to $N$. Thus under condition \eqref{trunc_M} both operators $L$ and $M$ preserve the space of polynomials of degree less or equal to $N$ and one can restrict the action of these operators to the $N+1$-dimensional space spanned by these polynomials. There are two natural bases for this space. The first one has the polynomials $\widehat{P}_n^{(\alpha,\beta)}(x)$ for $n=0,1,\ldots, N$ given in \eqref{Jac_hyp} as basis vectors; this is the eigenbasis of $L$. The second one has the functions $\psi_n(x)$ for $n=0,1,\ldots,N$ given by \eqref{psi_M} as basis vectors; this is the eigenbasis of $M$. In the parametrization \eqref{tau_nu}, the truncation condition \eqref{trunc_M} reads
\begin{equation*}
(\nu_2-N-1)(N+1+\nu_1 +\alpha+\beta) =0.
\end{equation*}      
It follows that one can take either $\nu_2=N+1$ or $\nu_1=-N-1-\alpha-\beta$. For definiteness, we choose to use the condition $\nu_2=N+1$. This leads to the following expression for the basis elements $\psi_n(x)$:
\begin{equation}
\label{psi_N_dim}
\psi_n(x) = x^N \: \widehat{P}_n^{(-\alpha-\beta-\nu_1-N-1,\beta)}(1/x) .
\end{equation}
It is easily seen that these are polynomials of degree $N-n$. Consider the expansion of the basis functions $\psi_n(x)$ given by \eqref{psi_N_dim} in terms of the basis functions $\widehat{P}_{k}^{(\alpha,\beta)}(x)$; one writes
\begin{equation*}
\psi_n(x) = \sum_{k=0}^N R_{nk} \widehat{P}_k^{(\alpha,\beta)}(x).
\end{equation*}
The expansion coefficients $R_{nk}$ can be obtained from the orthogonality relation \eqref{ort_Jac} of the basis functions $\widehat{P}_{n}^{(\alpha,\beta)}(x)$. One finds
\begin{align*}
 R_{nk} h_k = \int_0^1 \psi_n(x) \widehat{P}_k^{(\alpha,\beta)}(x)\; x^{\alpha} (1-x)^{\beta} dx,
\end{align*}
which can be considered as a finite-dimensional analog of the Fourier-Jacobi transform proposed by Koornwinder \cite{1985_Koornwinder_OPs&Appl_174}. In view of \eqref{Q_W_id}, this formula leads to an integral representation for the truncated Wilson (Racah) polynomials that arise in the expansion coefficients $R_{nk}$.  
\section{Tridiagonalization and the Racah--Wilson algebra}
In this section, the relation between the tridiagonalization of the hypergeometric operator, the Jacobi and Racah--Wilson algebras is discussed. It is shown that the operators $L$ and $M$ provide a realization of the Racah--Wilson algebra. The connection between this realization and the one arising from the Racah problem for the Lie algebra $\mathfrak{su}(1,1)$ is also established.
\subsection{The Racah--Wilson algebra}
The Racah--Wilson algebra \cite{1992_Granovskii&Lutzenko&Zhedanov_AnnPhys_217_1} is generated by the elements $A_1$, $A_2$ and their commutator $A_3$. These elements satisfy the commutation relations
\begin{align}
\label{W_def}
 \begin{aligned}
 {}
 [A_1,A_2]&=A_3,
 \\
[A_2,A_3]&=\alpha_1\{A_1,A_2\} + \alpha_2 A_2^2 + \gamma_1 A_1 + \delta A_2 + \epsilon_1,
\\
[A_3,A_1]&=\alpha_2\{A_1,A_2\} + \alpha_1 A_1^2 + \gamma_2 A_2 + \delta A_1+ \epsilon_2,
 \end{aligned}
 \end{align}
where $\alpha_1, \alpha_2, \gamma_1, \gamma_2, \delta,\epsilon_1,\epsilon_2$ are real structure constants. The algebra \eqref{W_def} has a Casimir operator which has the expression
\begin{multline*}
 R=\alpha_1\{A_1^2,A_2\} + \alpha_2
\{A_2^2,A_1\} + (\alpha_1^2+\gamma_1) A_1^2 + (\alpha_2^2+\gamma_2) A_2^2 +A_3^2
\\
+(\delta+\alpha_1 \alpha_2) \{A_1,A_2\} +(\delta \alpha_1 + 2\epsilon_1) A_1 +(\delta \alpha_2+2 \epsilon_2)A_2,
\end{multline*}
and which commutes with all generators $A_1, A_2,A_3$. Upon comparing the defining relations \eqref{W_def} of the Racah--Wilson algebra and those of the quadratic Jacobi algebra \eqref{QR3}, it is directly observed that the latter is a special case of the former that corresponds to taking $\alpha_1=\gamma_1=\epsilon_1=0$. The Racah--Wilson algebra is the algebraic structure associated to the Racah and Wilson polynomials, as it can be realized in terms of the operators involved in the bispectrality of these families of orthogonal polynomials. This realization, which can be found in \cite{2014_Genest&Vinet&Zhedanov_LettMathPhys_104_931} for example, is similar in spirit to the realization of the Jacobi algebra \eqref{QR3} in terms of the operators $L$ and $X$ exhibited in Section 2. For a review of this algebra and of some of its applications, one can consult \cite{2014_Genest&Vinet&Zhedanov_JPhysConfSer_512_012011}.
\subsection{Racah--Wilson algebra arising from tridiagonalization}
We now show how the Racah--Wilson algebra arises in the present context. Let $A_1(x)$ and $A_2(x)$ be defined as
\begin{align}
\label{Ultra}
 A_1(x)=L,\qquad A_2(x)=M,
\end{align}
where $L$ and $M$ are as in \eqref{hyp_op} and \eqref{M-Def}, respectively. Moreover, let $A_3(x)$ read
\begin{align*}
 A_3(x)=[A_1(x), A_2(x)].
\end{align*}
Then, by a direct calculation, it is verified that the operators $A_1(x)$, $A_2(x)$, $A_3(x)$ satisfy the defining relations \eqref{W_def} of the Racah--Wilson algebra with the structure constants taking the values
\begin{align*}
 \begin{aligned}
  &\alpha_1=-2,\qquad \alpha_2=2,\qquad \delta=\alpha(\alpha+\beta+1)+\beta-4\tau_0-2\tau_3,
  \\
  & \gamma_1=1-\alpha^2+4(\tau_0+\tau_2-\tau_2^2),\qquad \gamma_2=-(\alpha+\beta)(\alpha+\beta+2),
  \\
  & \epsilon_1=2\tau_0^2+(1-\alpha^2)\tau_3+2\tau_0\tau_3+2(\alpha+1)(\beta+1)(\tau_2^2-\tau_2)-(\alpha+1)(\alpha+\beta)\tau_0,
  \\
  &\epsilon_2= (\alpha+\beta)((\alpha+\beta+2)\tau_0 + (\alpha+1)\tau_3).
 \end{aligned}
\end{align*}
In this realization, the Casimir operator acts as a multiple of the identity
\begin{multline*} 
R(x)=2(\alpha+1)(\beta+1)(\alpha+\beta)(\tau_1^2 - \tau_1) + ((\alpha+\beta+1)^2-5) \tau_0^2 
\\
+2(\alpha+1)(\alpha+\beta) (\tau_0 \tau_3 + \tau_0 + \tau_3) - 4 \tau_0 \tau_3 +(\alpha^2-1)\tau_3^2.
\end{multline*}
In terms of the generators \eqref{Rea-1} of the quadratic Jacobi algebra \eqref{QR3}, the generators \eqref{Ultra} of the Racah--Wilson algebra have the expression
\begin{align*}
 A_1(x)=K_1,\quad A_2(x)=\tau_1 K_2\,K_1+\tau_2 K_1\,K_2+\tau_3 K_2+\tau_0.
\end{align*}
This provides an embedding of the Racah algebra into the Jacobi algebra in their standard realization on the space of polynomials.
\begin{remark}
In lieu of $M$, one can introduce an operator $J$ by applying the tridiagonalization procedure to the hypergeometric operator $L$ using instead the ``multiplication by $x-1$'' operator; this operator is defined by means of
\begin{equation}
 \label{J_def}
J =\rho_1(X-1)L+\rho_2L(X-1)+\rho_3(X-1)+\rho_0,
\end{equation}
where $\rho_1, \rho_2, \rho_3, \rho_0$ are arbitrary parameters such that $\rho_1+\rho_2 \ne 0$. Due to the symmetry of the hypergeometric operator with respect to the change of variable $x \to 1-x$, it is clear that the operator $J$ has the same properties as the operator $M$. In particular, one obtains the generic Wilson polynomials when expanding the eigenfunctions of the operator $J$ in a series of Jacobi polynomials. Furthermore, it is verified that if one defines the operators $B_1(x)=L$, and $B_2(x)=J$ as well as $B_3(x)=[B_1(x), B_2(x)]$, one obtains another realization of the Racah--Wilson algebra \eqref{W_def}.
\end{remark}
\subsection{Connection with the Racah problem}
There is an interesting interpretation of these results in terms of the Racah problem for the $\mathfrak{su}(1,1)$ algebra. Consider the $\mathfrak{su}(1,1)$ algebra with generators $S_0, S_{+}, S_{-}$ and commutation relations 
\begin{equation*}
[S_0, S_{\pm}]= \pm S_{\pm}, \quad [S_-, S_+]=2 S_0.
\end{equation*}
The Casimir operator
\begin{equation*}
C \equiv S_0^2 - S_0 - S_+S_- ,
\end{equation*}
commutes with all generators. On the irreducible representations of $\mathfrak{su}(1,1)$ that are relevant to the present paper, the Casimir operator takes the value $\sigma(\sigma-1)$.  For $\sigma>0$, the irreducible representations of the positive-discrete series are defined on the space spanned by the basis vectors $e_{n}$, $n=0,1,\ldots$, by the following actions:
\begin{equation*}
S_0 e_n = (n+\sigma) e_n,\quad  S_+ e_n = \gamma_{n+1} e_{n+1}, \quad S_- e_n = \gamma_n e_{n-1}, 
\end{equation*} 
where $\gamma_n = \sqrt{n(2 \sigma +n-1)}$. Consider the tensor product of three representations of the positive-discrete series with representation parameters $\sigma_1$, $\sigma_2$ and  $\sigma_3$. We use the notation
\begin{align*}
 S^{(ij)}_0=S^{(i)}_0+S^{(j)}_0,\qquad S^{(ij)}_{\pm}=S^{(i)}_{\pm}+S^{(j)}_{\pm},\qquad i,j=1,2,3,
\end{align*}
where the superscript $i$ in $S^{(i)}$ specifies on which representation space the generator acts.
The intermediate Casimir operators $C_{ik}$ are defined as
\begin{equation*}
C_{ik} = [S_0^{(ik)}]^2 - S_0^{(ik)} - S_+^{(ik)}S_-^{(ik)},\qquad\text{for}\quad (ik)\in \{(12),(23), (13)\},
\end{equation*}
and the total Casimir operator $C$ is written as
\begin{equation*}
C=[S_0^{(123)}]^2 - S_0^{(123)} - S_+^{(123)}S_-^{(123)},
\end{equation*}
where $S_{\pm}^{(ijk)}=S^{(i)}_{\pm}+S^{(j)}_{\pm}+S^{(k)}_{\pm}$ and similarly for $S_0^{(ijk)}$. By construction, the total Casimir operator commutes with all intermediate Casimir operators. 

Consider the decomposition of the three-fold tensor product representation in irreducible components. On each of these components, the total Casimir operator takes the value $\sigma_4(\sigma_4-1)$, where $\sigma_4=N+\sigma_1+\sigma_2+\sigma_3$ and $N$ is a non-negative integer. For a given $N$, there corresponds $N+1$ eigenvalues of the intermediate Casimir operators $C_{ij}$ which are of the form $\sigma_{ij}(\sigma_{ij}-1)$ where $\sigma_{ij}=k_{ij}+\sigma_i+\sigma_j$ and where $k_{ij}=0,1,\ldots,N$. Let $e_{k_{12}}^{(12)}$ and $e_{k_{23}}^{(23)}$ be eigenbases for $C_{12}$ and $C_{23}$, i.e.
\begin{equation*}
C_{12} e_{k_{12}}^{(12)} = \sigma_{12} (\sigma_{12} -1) e_{k_{12}}^{(12)}, \quad C_{23} e_{k_{23}}^{(23)} = \sigma_{23} (\sigma_{23} -1) e_{k_{23}}^{(23)},
\end{equation*}
where $ k_{12}, k_{23}=0,1,\dots,N$. The Racah problem consists in obtaining the expansion coefficients of the basis $e_{k_{23}}^{(23)}$ over the basis $e_{k_{12}}^{(12)}$
\begin{equation*}
e_{k_{23}}^{(23)} = \sum_{k_{12}=0}^N R_{k_{12}k_{23}} e_{k_{12}}^{(12)}.
\end{equation*}
For more details on the Racah problem for $\mathfrak{su}(1,1)$, the reader is referred to \cite{2003_VanderJeugt_OrthPoly&SpecFuns_25}.

A one-dimensional model for the Racah problem was constructed in \cite{2014_Genest&Vinet&Zhedanov_JPhysA_47_025203}  by separation of variables and dimensional reduction. It allows to express the three intermediate Casimir operators as second-order differential operators in one variable. One has
\begin{align}
\label{S}
\begin{aligned}
 C_{12}&= x^2(1-x)\, \partial_x^2 + x[(N-1-2 \sigma_1)x + 2(\sigma_1 + \sigma_2)]\partial_x
 \\
 &\qquad \qquad +2N \sigma_1 x +(\sigma_1+\sigma_2)(\sigma_1+\sigma_2-1),
 \\[.2cm]
 C_{23}&= x(x-1)\, \partial_x^2 + [2(1-N-\sigma_2 -\sigma_3)x + (N-1+2\sigma_3)]\partial_x 
 \\
 &\qquad \qquad +(N+\sigma_3+\sigma_2)(N+\sigma_3+\sigma_2-1),
 \\[.2cm]
 C_{31}&=x(x-1)^2 \partial_x^2 + (1-x)[(N-1-2\sigma_1)x + 1-N-2\sigma_3]\partial_x 
 \\
 &\qquad \qquad +2N \nu_1 (1-x) + (\sigma_3+\sigma_1)(\sigma_3+\sigma_1-1) .
 \end{aligned}
\end{align}
Any pair of operators $C_{ik}$ satisfy the Racah-Wilson algebra \cite{2014_Genest&Vinet&Zhedanov_JPhysA_47_025203}. Up to affine transformations, the operators $L,M$ and $J$ introduced in \eqref{hyp_op}, \eqref{M-Def} and \eqref{J_def}, coincide with the operators $C_{ik}$ of the one-dimensional model \eqref{S}. This provides an explanation for the fact that these operators satisfy the Racah--Wilson algebra. Moreover, it shows that the tridiagonalization procedure of the hypergeometric operator is intimately related to the Racah problem for $\mathfrak{su}(1,1)$.
\section{Degenerate case. The Hahn algebra}
In this section, we consider the degenerate case for which $\tau_1+\tau_2=0$ in the tridiagonalization procedure. It is shown that in this case the Hahn algebra and Hahn polynomials arise.

Let $\tau_1+\tau_2=0$ in \eqref{M-Def}. One can set $\tau_1=1$ and $\tau_2=-1$. With these conditions, the operator $M$ has the expression
\begin{equation}
 \label{M_deg}
M=2 x(x-1) \partial_x + (\alpha+\beta + \tau_3+2)x + \tau_0 -\alpha-1.
\end{equation}
It is observed that $M$ is now a first-order differential operator. Let 
\begin{equation}
\label{tau_3_N} 
\tau_3 = -2N-2-\alpha-\beta,
\end{equation}
with $N$ a non-negative integer. It is seen that both operators $L$ and $M$ preserve the $N+1$-dimensional space of polynomials of degree less or equal to $N$. It is already known that the eigenfunctions of $L$ are the polynomials \eqref{Jac_hyp}. The polynomial eigenfunctions $\psi_n(x)$ of the operator $M$ have the expression
\begin{equation}
\psi_n(x) = x^n (1-x)^{N-n}, \label{psi_Hahn} 
\end{equation}
and satisfy the eigenvalue equation
\begin{equation}
M\,\psi_n(x) = (\tau_0 -2n-\alpha-1)\,\psi_n(x),\label{M_psi_Hahn} 
\end{equation} 
for $n=0,1,\ldots, N$. By construction, the operator $M$ is 3-diagonal on the polynomial basis $\widehat{P}_n^{(\alpha,\beta)}(x)$: 
\begin{multline}
M \widehat{P}_n^{(\alpha,\beta)}(x) = (\lambda_n-\lambda_{n+1} + \tau_3) \widehat{P}_{n+1}^{(\alpha,\beta)}(x) 
\\+ (\tau_3 b_n+\tau_0) \widehat{P}_n^{(\alpha,\beta)}(x) + u_n(\lambda_n-\lambda_{n-1} + \tau_3) \widehat{P}_{n-1}^{(\alpha,\beta)}(x), \label{M_P_Hahn}
\end{multline}
where $\lambda_n$ is given by \eqref{lambda_Jac}. It is easy to show that the operator $L$ acts in a tridiagonal fashion on $\psi_n(x)$. Indeed, one has
\begin{multline}
L \psi_n(x) = (N-n)(N-n+\beta) \psi_{n+1}(x) 
\\
+
(2n^2 +(\alpha-\beta-2N)n - (\alpha+1)N) \psi_n(x)+n(n+\alpha)\psi_{n-1}.  \label{L_psi_Hahn}
\end{multline}
Applying the same method as in Section 2, one arrives at orthogonal polynomials $Q_n(\Lambda)$ satisfying the three-term recurrence relation \eqref{rec_Q} with coefficients
\begin{equation}
B_n =  \tau_3 b_n+\tau_0, \quad U_n = u_n (\tau_3 + \lambda_n -\lambda_{n-1})(\tau_3 + \lambda_{n-1} -\lambda_{n}). \label{BU_Hahn} \end{equation}
A direct calculation yields
\begin{multline*}
 B_n=\tau_0-N-1 -\frac{\alpha+\beta}{2}
 \\+
\frac{(\alpha^2-\beta^2)(2N+2+\alpha+\beta)}{4} \left( \frac{1}{\alpha+\beta+2n+2} - \frac{1}{\alpha+\beta+2n} \right) ,
\end{multline*}
and
\begin{equation*}
U_n = \frac{4n(n+\alpha)(n+\beta)(n+\alpha+\beta)(n+N+1+\alpha+\beta)(N+1-n)}{(2n+\alpha+\beta-1)(2n+\alpha+\beta)^2(2n+\alpha+\beta+1)}.
\end{equation*}
It is seen that these recurrence coefficients correspond to the generic Hahn polynomials $Q_n(x;\alpha,\beta;N)$; see \cite{2010_Koekoek_&Lesky&Swarttouw} for the definition. 

The operators $L$ and $M$ defined as in \eqref{hyp_op} and \eqref{M_deg} also generate a quadratic algebra. Let $Z$ be the operator defined by 
\begin{equation*}
Z= [L,M] = 2x(1-x)(2x-1) \partial_x^2 - 2V_1(x) \partial_x - V_0(x),  
\end{equation*}
where
\begin{align*}
V_1(x)&= (\tau_3 + 2 \alpha +6 + 2 \beta)x^2 -(\tau_3 +6 + 3 \alpha +\beta)x + \alpha+1, 
\\
V_0(x)&= (\alpha+\beta+2 +\tau_3)((\alpha+\beta+2)x-\alpha-1).
\end{align*}
It is verified that $L$, $M$ and $Z$ satisfy the commutation relations
\begin{align*}
[M,Z]&=2\,M^2 -2(\tau_3 + 2 \tau_0) M - 4\,L + 2\tau_0(\tau_0+\tau_3) + 2 (\alpha+1)(\beta+1),
\\[.2cm]
[Z,L]&=2\{L,M\} -2(\tau_3 + 2 \tau_0) L -(\alpha+\beta)(\alpha+\beta+2)M 
\\
& \qquad \qquad \qquad\qquad +(\alpha+\beta)(\tau_0(\alpha+\beta+2) + \tau_3(\alpha+1)). 
\end{align*}
These relations correspond to the quadratic Hahn algebra \cite{1992_Granovskii&Lutzenko&Zhedanov_AnnPhys_217_1}, which is a special case of the Racah--Wilson algebra \eqref{W_def}.
\section{Conclusion}
In this paper, we have explained from an algebraic standpoint how the tridiagonalization of the hypergeometric operator leads to the Wilson polynomials; that is we have shown how the operators involved generate the Racah--Wilson algebra.

It would be of interest to enlarge the scope of our study to cover $q$-orthogonal polynomials and their $q\rightarrow -1$ limits. In this context, it would be desirable to arrive at the Bannai--Ito algebra by the tridiagonalization approach and to determine how the Bannai--Ito polynomials can be obtained from lower families of orthogonal polynomials of the $-1$ hierarchy.
\section*{Acknowledgments}
\noindent
The authors would thank two referees for valuable remarks. The authors would also like to thank E. Koelink for stimulating discussions and for his suggestions on the manuscript. V. X. G. holds a fellowship from the Natural Science and Engineering Research Council of Canada (NSERC). The research of M. E. H. I. is supported by the DSFP of King Saud University.  The research of L. V. is supported in part by NSERC. A. Z. wishes to thank the Centre de recherches math\'ematiques (CRM) for its hospitality.

\end{document}